\documentstyle[11pt,twoside]{article}
\include{graphicx}
\title{{\bf The Stokes phenomenon and the Lerch zeta function}}
\author{\sc R. B.\ Paris\footnote{E-mail address:\ \ {\tt r.paris@abertay.ac.uk}}\\
{\em Division of Computing and Mathematics}, \\
{\em University of Abertay Dundee, Dundee DD1 1HG, UK}\\
}
\begin{document}
\def\f#1#2{\mbox{${\textstyle \frac{#1}{#2}}$}}
\def\dfrac#1#2{\displaystyle{\frac{#1}{#2}}}
\def\boldal{\mbox{\boldmath $\alpha$}}
\newcommand{\bee}{\begin{equation}}
\newcommand{\ee}{\end{equation}}
\newcommand{\lam}{\lambda}
\newcommand{\ka}{\kappa}
\newcommand{\al}{\alpha}
\newcommand{\om}{\omega}
\newcommand{\Om}{\Omega}
\newcommand{\fr}{\frac{1}{2}}
\newcommand{\fs}{\f{1}{2}}
\newcommand{\g}{\Gamma}
\newcommand{\br}{\biggr}
\newcommand{\bl}{\biggl}
\newcommand{\ra}{\rightarrow}
\newcommand{\gtwid}{\raisebox{-.8ex}{\mbox{$\stackrel{\textstyle >}{\sim}$}}}
\newcommand{\ltwid}{\raisebox{-.8ex}{\mbox{$\stackrel{\textstyle <}{\sim}$}}}
\renewcommand{\topfraction}{0.9}
\renewcommand{\bottomfraction}{0.9}
\renewcommand{\textfraction}{0.05}
\newcommand{\mcol}{\multicolumn}
\date{}
\maketitle
\begin{abstract}
We examine the exponentially improved asymptotic expansion of the Lerch zeta function $L(\lambda,a,s)=\sum_{n=1}^\infty \exp (2\pi ni\lambda)/(n+a)^s$ for large complex values of $a$, with $\lambda$ and $s$ regarded as parameters. It is shown that an infinite number of subdominant exponential terms switch on across the Stokes lines $\arg\,a=\pm\fs\pi$. In addition, it is found that the transition across the upper and lower imaginary $a$-axes is associated, in general, with unequal scales. Numerical calculations are presented to confirm the theoretical predictions.
\vspace{0.4cm}

\vspace{0.1cm}

\end{abstract}

\vspace{0.3cm}
\begin{center}
{\bf 1.\ Introduction}
\end{center}
\setcounter{section}{1}
\setcounter{equation}{0}
\renewcommand{\theequation}{\arabic{section}.\arabic{equation}}
A principal result of asymptotic analysis in the last quarter century  has been the interpretation of the Stokes phenomenon as the smooth appearance of an exponentially small term in compound asymptotic expansions across certain rays in the complex plane known as Stokes lines. For a wide class of functions, particularly those satisfying second-order ordinary differential equations, the functional form of the coefficient multiplying such a subdominant exponential (a Stokes multiplier) is found to possess a universal structure represented to leading order by an error function,
whose argument is an appropriate variable describing the transition across the Stokes line \cite{B1}.

A function not satisfying a differential equation and which does not share this simple property is the logarithm of the gamma function. In \cite{PW}, Paris and Wood obtained the exponentially improved expansion of $\log\,\g(z)$ and showed that it involved not one but an {\it infinite number} of subdominant exponentials $e^{\pm 2\pi ikz}$ ($k=1, 2, \ldots$). These exponentials are maximally subdominant on the Stokes lines $\arg\,z=\pm\fs\pi$, respectively, and steadily grow in magnitude in $|\arg\,z|>\fs\pi$ to eventually combine to generate the poles of $\log\,\g(z)$ on the negative $z$-axis. These authors demonstrated that the Stokes multipliers associated with the leading exponentials 
(corresponding to $k=1$) undergo a smooth transition approximately described by an error function in the neighbourhood of $\arg\,z=\pm\fs\pi$.
Subsequently, Berry \cite{B2} showed, by a sequence of increasingly delicate subtractions of optimally truncated asymptotic series, that all the subdominant exponentials switch on smoothly across the Stokes lines with the multiplier given approximately by
\bee\label{e11}
\hspace{4cm}\fs\pm\fs\,\mbox{erf}\,[(\theta\mp\fs\pi)\sqrt{\pi k|z|}\,], \qquad(k=1, 2, \ldots )
\ee
in the neighbourhood of $\theta=\arg\,z=\pm\fs\pi$, respectively; see also \cite[\S 6.4]{PK} for a detailed summary.

An analogous refinement in the large-$a$ asymptotics of the Hurwitz zeta function 
\[\zeta(s,a)=\sum_{n=0}^\infty (n+a)^{-s}\qquad (\Re (s)>1;\ a\neq 0, -1, -2, \ldots)\]
was considered in \cite{P}. Across the Stokes lines $\arg\,a=\pm\fs\pi$, there is a similar appearance of an infinite number of subdominant exponentials $e^{\pm 2\pi ika}$ ($k=1, 2, \ldots $), each exponential being associated with its own Stokes multiplier. For large $|a|$, the Stokes multipliers associated with these exponentials   
also undergo a smooth, but rapid, transition in the neighbourhood of $\arg\,a=\pm\fs\pi$ given approximately by 
(\ref{e11}) with $z$ replaced by $a$

In \cite{P1}, the periodic zeta function $F(\lambda,s)$, given by \cite[\S 25.13]{DLMF}
\bee\label{e12}
F(\lambda,s)=\sum_{n=1}^\infty \frac{e^{2\pi ni\lambda}}{n^s}\qquad (\Re (s)>0,\ 0<\lambda<1;\  \Re (s)>1,\ \lambda\in{\bf N}),
\ee
was discussed for complex values of the parameter $\lambda$ in the upper half-plane $0<\arg\,\lambda<\pi$.
This function can be expressed in terms of the Hurwitz zeta function and, accordingly, its exponentially improved large-$\lambda$ expansion also consists of an infinite number of subdominant exponentials $e^{2\pi ik\lambda}$ ($k=\pm 1, \pm 2,\ldots$). In the neighbourhood of the positive imaginary $\lambda$-axis, it is found that the exponentials with $k\geq 1$ undergo a double Stokes phenomenon, since constituent parts of $F(\lambda,s)$ are associated with two parallel Stokes lines at unit distance apart.

In this paper we consider the Lerch zeta function defined for $\Re (s)>1$, $0<\lambda\leq 1$ and complex $a$ by the series \cite[Eq.~(25.14.1)]{DLMF}
\bee\label{e13}
L(\lambda,a,s):=\sum_{n=0}^\infty \frac{e^{2\pi ni\lambda}}{(n+a)^s}\qquad (a\neq 0, -1, -2, \ldots)
\ee
and elsewhere by analytic continuation. When $\lambda=1$, the Lerch function reduces to $\zeta(s,a)$ and when $a=1$ we have
\bee\label{e14}
L(\lambda,1,s)=e^{-2\pi i\lambda} F(\lambda,s).
\ee
We shall find in the large-$a$ asymptotics of $L(\lambda,a,s)$ that there is a similar
appearance of an infinite number of subdominant exponential terms in the neighbourhood of the Stokes lines $\arg\,a=\pm\fs\pi$. Each of these exponentials is associated with its own Stokes multiplier, which undergoes a smooth, but rapid, transition in the vicinity of these rays. However, unlike the situation present with $\zeta(s,a)$, it will be found that the transition across the Stokes lines in the upper and lower half-planes is, in general, associated with unequal scales.  

We use a Mellin-Barnes integral definition of $L(\lambda,a,s)$ to first determine its large-$a$ Poincar\'e expansion
and then its exponentially improved expansion as $|a|\ra\infty$ in the sector $|\arg\,a|<\pi$. The procedure we adopt is similar to that employed in \cite{P, PW}.

\vspace{0.6cm}

\begin{center}
{\bf 2.\ An integral representation}
\end{center}
\setcounter{section}{2}
\setcounter{equation}{0}
\renewcommand{\theequation}{\arabic{section}.\arabic{equation}}
Let $a$ and $s$ be complex variables with $|\arg\,a|<\pi$ and $\lambda$ a real variable satisfying $0<\lambda\leq 1$. When $\Re (s)>1$, we can write the Lerch zeta function defined in (\ref{e13}) as
\[L(\lambda,a,s)=a^{-s}+a^{-s}\sum_{n=1}^\infty e^{2\pi i\lambda} \bl(1+\frac{n}{a}\br)^{-s}.\]
Making use of the representation (see, for example, \cite[p.~91]{PK})
\[\frac{\g(s)}{(1+x)^s}=\frac{1}{2\pi i}\int_{c-\infty i}^{c+\infty i}  \g(u) \g(s-u)x^{-u}du\qquad (|\arg\,x|<\pi),\]
where\footnote{This integral representation holds for the wider $c$-interval given by $0<c<\Re (s)$.} $1<c<\Re (s)$, we obtain after an interchange in the order of summation and integration
\bee\label{e21}
L(\lambda,a,s)=a^{-s}+\frac{a^{-s}}{2\pi i\,\g(s)}\int_{c-\infty i}^{c+\infty i} \g(u) \g(s-u) F(\lambda,u)\,a^u du,
\ee
where $F(\lambda,u)$ is the periodic zeta function defined in (\ref{e12}). When $\lambda=1$, $F(1,u)$ becomes the Riemann zeta function $\zeta(u)$ and the representation (\ref{e21}) reduces to that given in \cite{P} for the Hurwitz zeta function $L(1,a,s)=\zeta(s,a)$.

Displacement of the integration path in (\ref{e21}) to the left over the simple pole at $u=0$ (and, in the case $\lambda=1$, at $u=1$) then yields
\bee\label{e22}
L(\lambda,a,s)=
\frac{\epsilon(\lambda) a^{1-s}}{s-1}
+a^{-s}\{1+F(\lambda,0)\}
+\frac{Z(\lambda,a,s)}{\g(s)},
\ee
where $\epsilon(\lambda)=0$ or 1 according as $0<\lambda<1$ or $\lambda=1$ and, with the change of variable $u\ra -u$,
\bee\label{e23}
Z(\lambda,a,s)=\frac{a^{-s}}{2\pi i}\int_{c-\infty i}^{c+\infty i} \g(-u) \g(u+s) F(\lambda,-u)\,a^{-u} du\qquad (0<c<1).
\ee

The result in (\ref{e22}) and (\ref{e23}) has been derived assuming that $\Re (s)>1$; but this restriction can be relaxed to allow for $\Re (s)\leq 1$ by suitable indentation of the integration path to lie to the right of all the poles of $\g(s+u)$ (provided $s\neq -1, -2, \ldots $). The representation (\ref{e22}) is similar to, but not identical with, that given in \cite{K}.
\vspace{0.3cm}

\noindent{\it 2.1.\ The Poincar\'e asymptotic expansion}
\vspace{0.1cm}

\noindent
The large-$|a|$ asymptotic expansion of $L(\lambda,a,s)$ can be obtained by further displacement of the integration path over the poles of $\g(-u)$ at $u=1, 2, \ldots , K-1$, where $K$ denotes an arbitrary positive integer. This  yields, when $0<\lambda<1$,
\bee\label{e24}
L(\lambda,a,s)=a^{-s}+\sum_{k=0}^{K-1}\frac{(-)^k}{k!} (s)_k F(\lambda,-k)\,a^{-s-k}+{\cal R}_K(\lambda,a,s),
\ee
where $(a)_k=\g(a+k)/\g(a)$ is Pochhammer's symbol. In the Appendix A it is shown that the remainder ${\cal R}_K(\lambda,a,s)=O(a^{-s-K})$ as $|a|\ra\infty$ in $|\arg\,a|<\pi$. 
This expansion agrees with that given by Ferreira and L\'opez \cite[Thm 1]{FL}, who expressed their coefficients in terms of the polylogarithm function $\mbox{Li}_{-n}(e^{2\pi i\lambda})=F(\lambda,-n)$.

In \cite{A}, it is shown that $F(\lambda,-k)$ is expressible in terms of a certain kind of generalised Bernoulli polynomials ${\tilde B}_n(\alpha,\beta)$ defined by the generating function
\[\frac{te^{\alpha t}}{1-\beta e^t}=\sum_{n=0}^\infty \frac{{\tilde B}_n(\alpha,\beta)}{n!}\,t^n\]
in a neighbourhood of $t=0$. From (\ref{e14}) and \cite[p.~164]{A}, we have for non-negative integer values of $k$ 
\bee\label{e25a}
F(\lambda,-k)=e^{2\pi i\lambda}L(\lambda,1,-k)=\frac{e^{2\pi i\lambda}}{k+1}\,{\tilde B}_{k+1}(1,e^{2\pi i\lambda}),
\ee
where ${\tilde B}_n(1,x)$ may be expressed in the form 
\bee\label{e25}
{\tilde B}_n(1,x)=\frac{n P_n(x)}{(1-x)^n},
\ee
with $P_1(x)=1$ and for $n\geq 2$ the $P_n(x)$ are polynomials of degree $n-2$. From \cite[Eq.~(3.7)]{A}, we have that
\[P_n(x)=\sum_{r=1}^{n-1}  r! {\bf S}_{n-1}^{(r)}\,x^{r-1} (1-x)^{n-r-1}\qquad(n\geq 2)\]
where ${\bf S_n^{(r)}}$ are the Stirling numbers of the second kind. The first few $P_n(x)$ are consequently:
\begin{eqnarray*}
P_2(x)&=&1,\quad P_3(x)=1+x,\quad P_4(x)=1+4x+x^2,\\
P_5(x)&=&1+11x+11x^2+x^3,\quad P_6(x)=1+26x+66x^2+26x^3+x^4,\\
P_7(x)&=&1+57x+302x^2+302x^3+57x^4+x^5,\\
P_8(x)&=&1+120x+1191x^2+2416x^3+1191x^4+120x^5+x^6, \ldots\,.
\end{eqnarray*}

Then we have the expansion given in the following theorem.
\newtheorem{theorem}{Theorem}
\begin{theorem}$\!\!\!.$ Let $s$ $(\neq 0, -1 -2, \ldots)$ be a complex variable and $K$ be an arbitrary positive integer. Then, when $0<\lambda<1$, we have 
\bee\label{e26}
L(\lambda,a,s)=\frac{a^{-s}}{1-e^{2\pi i\lambda}}+e^{2\pi i\lambda}\sum_{k=1}^{K-1}\frac{(-)^k (s)_k}{k!}\,\frac{P_{k+1}(e^{2\pi i\lambda})\,a^{-s-k}}{(1-e^{2\pi i\lambda})^{k+1}}+O(a^{-s-K})
\ee
as $|a|\ra\infty$ in the sector $|\arg\,a|<\pi$.
\end{theorem}

It is worth remarking that the first few coefficients of the expansion (\ref{e26}) can be expressed in an alternative trigonometric form to yield
\[L(\lambda,a,s)\sim \frac{a^{-s}}{2\sin \pi\lambda}\bl\{ie^{-\pi i\lambda}+\frac{s}{2a\sin \pi\lambda}-\frac{is(s+1) \cos \pi\lambda}{4a^2\sin^2\! \pi\lambda}-\frac{s(s+1)(s+2)(1+2\cos^2\!\pi\lambda)}{23a^3\sin^3\! \pi\lambda}\]
\[+\frac{is(s+1)(s+2)(s+3)(2+\cos^2\!\pi\lambda)\cos \pi\lambda}{48a^4\sin^4\!\pi\lambda}+\cdots\br\}.\]

When $\lambda=1$, we have $F(\lambda,-k)=\zeta(-k)$ in (\ref{e24}). Then, from (\ref{e22}) and the fact that  $\zeta(0)=-\fs$, $\zeta(-2k)=0$ and $\zeta(-2k+1)= -B_{2k}/(2k)$ ($k=1, 2, \ldots$), where $B_{2k}$ denote the even-order Bernoulli numbers, we recover the well-known asymptotic expansion \cite[p.~25]{MO}
\[L(1,a,s)\equiv \zeta(s,a)\sim \frac{1}{2}a^{-s}+\frac{a^{1-s}}{s-1}+\frac{1}{\g(s)} \sum_{k=1}^\infty \frac{B_{2k}}{(2k)!}\,\frac{\g(2k+s-1)}{a^{2k+s-1}}\]
valid as $|a|\ra\infty$ in $|\arg\,a|<\pi$.

\vspace{0.6cm}

\begin{center}
{\bf 3.\ The exponentially improved expansion}
\end{center}
\setcounter{section}{3}
\setcounter{equation}{0}
\renewcommand{\theequation}{\arabic{section}.\arabic{equation}}
Let $0<\lambda\leq 1$ and define $\lambda'=1-\lambda$. To determine the exponentially improved expansion of $L(\lambda,a,s)$ for large $|a|$ in $|\arg\,a|<\pi$ and, in particular, the behaviour in the neighbourhood of the rays $\arg\,a=\pm\fs\pi$, we start with the representation in (\ref{e22}), namely 
\bee\label{e31}
L(\lambda,a,s)=\frac{\epsilon(\lambda)a^{1-s}}{s-1}+a^{-s}\{1+F(\lambda,0)\}+\frac{Z(\lambda,a,s)}{\g(s)},
\ee
where from (\ref{e25a}) and (\ref{e25}), we have $F(\lambda,0)=e^{2\pi i\lambda}/(1-e^{2\pi i\lambda})$ when $0<\lambda<1$ and $F(1,0)=-\fs$. 
The function $Z(\lambda,a,s)$ is defined by the integral in (\ref{e23}).
\vspace{0.3cm}

\noindent{\it 3.1.\ An expansion for $Z(\lambda,a,s)$ }
\vspace{0.1cm}

\noindent
The functional equation for the periodic zeta function $F(\lambda,s)$ can be obtained from the analogous result for $L(\lambda,a,s)$ given in \cite[pp.~26, 29]{E} and takes the form
\bee\label{e30}
F(\lambda,s)=\frac{\g(1-s)}{(2\pi)^s}\bl\{e^{\fr\pi i(1-s)}\sum_{k=0}^\infty(k+\lambda)^{s-1}+e^{-\fr\pi i(1-s)}\sum_{k=0}^\infty{}\!'(k+\lambda')^{s-1}\br\}\qquad (\Re (s)<0),
\ee
where the prime on the second summation sign indicates that the term corresponding to $k=0$ is to be omitted if $\lambda=1$.
Substitution of this last result into (\ref{e23}) yields
\[Z(\lambda,a,s)=\frac{a^{-s}}{2\pi i}\int_{c-\infty i}^{c+\infty i} \g(-u) \g(u+s) F(\lambda,-u)\,a^{-u} du\hspace{3cm}\]
\[
\hspace{0.8cm}=-\frac{a^{-s}}{4\pi} \int_{c-\infty i}^{c+\infty i} \frac{\g(u+s)}{\sin \pi u} (2\pi a)^{-u}\hspace{4cm}\]
\[\hspace{4cm}\times\bl\{e^{\fr\pi iu} \sum_{k=0}^\infty (k+\lambda)^{-1-u}-e^{-\fr\pi iu} \sum_{k=0}^\infty{}\!'(k+\lambda')^{-1-u}\br\} du,\]
where $0<c<1$.
We define the sum
\bee\label{e31a}
E(\lambda,s;z):=\sum_{k=0}^\infty (k+\lambda)^{s-1} J_k(\lambda,s;z),
\ee
where
\bee\label{e31b}
J_k(\lambda,s;z):=\frac{1}{4\pi}
\int_{c-\infty i}^{c+\infty i}\frac{\g(u+s)}{\sin \pi u}\,z^{-u-s}du\qquad (0<c<1,\ |\arg\,z|<\f{3}{2}\pi).
\ee
Then we obtain
\bee\label{e32}
Z(\lambda,a,s)=(2\pi)^s \{e^{\fr\pi is} E(\lambda',s;iX')-e^{-\fr\pi is}
E(\lambda,s;-iX)\},
\ee
where 
\bee\label{e32d}
X:=2\pi a(k+\lambda),\quad X':=2\pi a(k+\lambda')
\ee
and in the sum $E(\lambda',s;iX')$ the term corresponding to $k=0$ is understood to be omitted if $\lambda=1$.

We now displace the (possibly indented when $\Re (s)\leq 0$) integration path for $J_k(\lambda,s;-iX)$ in (\ref{e31b}) to the right over the poles situated at $u=1, 2, \ldots , N_k-1$, where the $\{N_k\}$ ($k\geq 1$) denote (for the moment) an arbitrary set of positive integers. This produces
\bee\label{e31c}
J_k(\lambda,s;-iX)=\frac{1}{2\pi i}\sum_{r=1}^{N_k-1} (-)^r \g(r+s)(-iX)^{-r-s}+R_k(\lambda,a;N_k),
\ee
where the remainder term is given by
\begin{eqnarray}
R_k(\lambda,a;N_k)&=&\frac{1}{4\pi}\int_{-c'+N_k-\infty i}^{-c'+N_k+\infty i} \frac{\g(u+s)}{\sin \pi u}\,(-iX)^{-u-s}du\nonumber\\
&=&\frac{(-)^{N_k}}{4\pi}\int_{-c'-\infty i}^{-c'+\infty i} \frac{\g(u+\nu_k)}{\sin \pi u}\,(-iX)^{-u-\nu_k}du,\label{e33}
\end{eqnarray}
with $0<c'<1$. In the last integral, we have replaced the integration variable $u$ by $u+N_k$ and have set $\nu_k:=N_k+s$. It is also tacitly assumed that $\Re (\nu_k)-c'>0$ so that the integration path in $R_k(\lambda,a;N_k)$ is not indented when $\Re (s)\leq 0$.

The integrals in (\ref{e33}) can be identified in terms of the so-called terminant function (see \cite[p.~67]{DLMF}), which is a multiple of the incomplete gamma function $\g(a,z)$ (or, equivalently, the exponential integral). Following the notation employed in \cite[p.~243]{PK}, we denote the terminant function by $T_\nu(z)$, where
\begin{eqnarray}
T_\nu(z)&=&e^{\pi i\nu}\frac{\g(\nu)}{2\pi i}\,\g(1-\nu,z)\nonumber\\
&=&\frac{e^{-z}}{4\pi}\int_{-c-\infty i}^{-c+\infty i} \frac{\g(u+\nu)}{\sin \pi u}\,z^{-u-\nu}du \qquad (0<c<1,\ |\arg\,z|<\f{3}{2}\pi),\label{e32a}
\end{eqnarray}
where, provided $\nu\neq 0, -1, -2, \ldots\,$, the integration path lies to the right of all the poles of $\g(u+\nu)$; see \cite[p.~178]{DLMF}. It then follows from (\ref{e33}) that
\bee\label{e35a}
R_k(\lambda,a;N_k)=e^{-iX-\pi is}\,T_{\nu_k}(-iX)\qquad (|\arg\,a|<\pi).
\ee

Proceeding in the same manner for the integral $J_k(\lambda',s;iX')$, with the set of integers $\{N_k\}$ replaced by the set $\{N_k'\}$ and the parameters $\nu_k$ by $\nu_k':=N_k'+s$, we find 
\bee\label{e31d}
J_k(\lambda',s;iX')=\frac{1}{2\pi i}\sum_{r=1}^{N_k'-1} (-)^r \g(r+s)(iX')^{-r-s}+R_k'(\lambda',a;N_k'),
\ee
where
\bee\label{e35b}
R_k'(\lambda',a;N_k')=e^{iX'-\pi is}\,T_{\nu_k'}(iX')\qquad (|\arg\,a|<\pi).
\ee
Then, from (\ref{e32}), we finally obtain
\[\frac{Z(\lambda,a,s)}{(2\pi)^s}=\sum_{k=0}^\infty (k+\lambda)^{s-1}\bl\{\frac{i}{2\pi }\sum_{r=1}^{N_k-1}\frac{(-i)^r \g(r+s)}{X^{r+s}}-e^{-\fr\pi is} R_k(\lambda,a;N_k)\br\}\]
\bee\label{e35}
\hspace{0.5cm}-\sum_{k=0}^\infty{}\!' (k+\lambda')^{s-1}\bl\{\frac{i}{2\pi}\sum_{r=1}^{N_k'-1}\frac{i^r \g(r+s)}{X'^{\,r+s}}-e^{\fr\pi is} R_k'(\lambda',a;N_k')\br\}
\ee
valid in $|\arg\,a|<\pi$.

When $\lambda=1$, we have $X=X'=2\pi ika$ and $N_k=N_k'$. The right-hand side of (\ref{e35}) then reduces to
\[\sum_{k=1}^\infty k^{s-1}\bl\{\frac{1}{\pi}\sum_{r=0}^{n_k-1} (-)^r\frac{\g(2r+s+1)}{X^{2r+s+1}}+{\cal R}_k(a;n_k)\br\},\]
where $n_k$ denotes an arbitrary positive integer and 
\[{\cal R}_k(a;n_k)=e^{-\pi is}\{e^{iX+\fr\pi is} T_{\nu_k}(iX)-e^{-iX-\fr\pi is} T_{\nu_k}(-iX)\},\quad \nu_k:=2n_k+s\]
as found in \cite[Eq.~(2.5), (2.6)]{P} for the Hurwitz zeta function $\zeta(s,a)$.

Assuming that $N_k<N_{k+1}$, $N_k'<N_{k+1}'$ ($k=0, 1, 2, \ldots$), we can write the double sum involving $\lambda$ in
(\ref{e35}) in the form  \cite{B2}, \cite[\S 6.4.3]{PK}
\begin{eqnarray}
\sum_{k=0}^\infty(k+\lambda)^{s-1}\sum_{r=0}^{N_k-1}\frac{(-i)^r \g(r+s)}{X^{r+s}}&=&\sum_{k=0}^\infty\sum_{r=0}^{N_k-1} \frac{(-i)^r \g(r+s)}{(2\pi a)^{r+s} (k+\lambda)^{r+1}}\hspace{3.4cm}\nonumber\\
&=&\sum_{m=0}^\infty \sum_{r=N_{m-1}}^{N_m-1} \frac{(-i)^r \g(r+s)}{(2\pi a)^{r+s}}\,\zeta(r+1,m+\lambda),\label{e37g}
\end{eqnarray}
where the sum over $k$ has been evaluated in terms of the Hurwitz zeta function and $N_{-1}=1$; see Appendix B for details. A similar rearrangement applies to the other double sum in (\ref{e35}) involving $\lambda'$ (with $N_{-1}'=1$).

If we define 
\[H_m(a;\lambda,\lambda'):=\sum_{r=N_{m-1}}^{N_m-1} \frac{(-i)^r \g(r+s)}{(2\pi a)^{r+s}}\,\zeta(r+1,m+\lambda)\hspace{4cm}\]
\bee\label{e36a}
\hspace{4cm}-\sum_{r=N_{m-1}'}^{N_m'-1} \frac{i^r \g(r+s)}{(2\pi a)^{r+s}}\,\zeta(r+1,m+\lambda'),
\ee
then, from (\ref{e31}), (\ref{e35a}), (\ref{e35b}) and (\ref{e35}), we finally obtain
\begin{theorem}$\!\!\!.$
Let $\lambda'=1-\lambda$, $\epsilon(\lambda)=0$ or 1 according  as $0<\lambda<1$ or $\lambda=1$. Let the truncation indices $N_k$ and $N_k'$ be increasing sets of positive integers with $\nu_k=N_k+s$, $\nu_k'=N_k'+s$. Then, for $0<\lambda\leq 1$, we have the expansion of $L(\lambda,a,s)$ given by
\[L(\lambda,a,s)=\frac{\epsilon(\lambda)a^{1-s}}{s-1}+a^{-s}\{1+F(\lambda,0)\}
+\frac{(2\pi)^s}{\g(s)}\sum_{m=0}^\infty{}\!' \bl\{\frac{i}{2\pi}\,H_m(a;\lambda,\lambda')\hspace{2cm}\]
\bee\label{e37c}
\hspace{3cm}-\frac{e^{-\fr\pi is}}{(m+\lambda)^{1-s}}\,R_m(\lambda,a;N_m)
+\frac{e^{\fr\pi is}}{(m+\lambda')^{1-s}}\,R_m'(\lambda',a;N_m')\br\}
\ee
valid in $|\arg\,a|<\pi$, where 
\[R_m(\lambda,a;N_m)=e^{-iX-\pi is} T_{\nu_m}(-iX), \quad R_m'(\lambda',a;N_m')=e^{iX'-\pi is} T_{\nu_m'}(iX')\]
and $X$, $X'$ are defined in (\ref{e32d}). The prime on the summation sign indicates that the terms corresponding to $m=0$ in the sums involving $\lambda'$ are to be omitted when $\lambda=1$.
\end{theorem}

\vspace{0.3cm}

\noindent{\it 3.2.\ The optimally truncated expansion and the Stokes multipliers}
\vspace{0.1cm}

\noindent
An important feature of (\ref{e35}) is that the Poincar\'e expansion in (\ref{e26}) has been decomposed into two $k$-sequences of component asymptotic series with scales $2\pi a(k+\lambda)$ and $2\pi a(k+\lambda')$, each associated with its own {\it arbitrary} truncation index $N_k$ and $N_k'$ and remainder terms $R_k(\lambda,a;N_k)$ and $R_k'(\lambda',a;N_k')$. From the large-argument asymptotics of the incomplete gamma function \cite[p.~179]{DLMF} 
\[\g(a,z)\sim z^{a-1}e^{-z} \qquad (|z|\ra\infty,\ |\arg\,z|<\f{3}{2}\pi)\]
and the first equation in (\ref{e32a}),
the sums involving the remainders are absolutely convergent, since the decay of the late terms is controlled by 
$k^{-N_k-1}$ and $k^{-N_k'-1}$. It then follows that the result in (\ref{e35}) is exact and that no further expansion process is required. 

The infinite sequences of exponentials $e^{2\pi ia(k+\lambda)}$ and $e^{-2\pi ia(k+\lambda')}$ ($k=0, 1, 2, \ldots$) are seen to emerge from the remainders in (\ref{e35}), or (\ref{e37c}), with the terminant functions $T_{\nu_k}(-iX)$
and $T_{\nu_k'}(iX')$, respectively, as coefficients. These exponentials are maximally subdominant on the negative and positive imaginary axes, respectively and steadily increase in magnitude as one approaches the negative real $a$-axis where they eventually combine to generate the singularities of $L(\lambda,a,s)$ at negative integer values of $a$.

If the truncation indices $N_k$ and $N_k'$ are now chosen to correspond to the optimal truncation values (i.e., truncation at or near the least term in the corresponding inner series over $r$ in (\ref{e35})), then it is easily shown that
\bee\label{e36}
N_k\simeq 2\pi (k+\lambda)|a|,\qquad N_k'\simeq 2\pi (k+\lambda')|a|.
\ee
In this case, $\nu_k=|X|+O(1)$, $\nu_k'=|X'|+O(1)$ and we see that the order and the argument of each terminant function appearing in (\ref{e35}) are approximately equal in the limit $|a|\ra\infty$. When $|\nu|\sim |z|\gg 1$, the function $T_\nu(z)$ possesses the asymptotic behaviour \cite{O}, \cite[\S 6.2.6]{PK}
\bee\label{e37}
T_\nu(z)\sim\left\{\begin{array}{lc}\!\! \dfrac{-ie^{(\pi-\phi)i\nu}}{1+e^{-i\phi}}\,\frac{e^{-z-|z|}}{\sqrt{2\pi |z|}}\{1+O(z^{-1})\}& -\pi+\delta\leq\phi\leq\pi-\delta \\
\\
\!\!\fs+\fs \mbox{erf}\,[c(\phi)(\fs|z|)^\fr]+O(z^{-\fr}e^{-\fr|z|c^2(\phi)}),& \delta\leq\phi\leq2\pi-\delta\end{array}\right.
\ee
where $\phi=\arg\,z$, $\delta$ denotes an arbitrarily small positive quantity and $c(\phi)$ is defined implicitly by
\[\fs c^2(\phi)=1+i(\phi-\pi)-e^{i(\phi-\pi)}\]
with the branch for $c(\phi)$ chosen so that $c(\phi)\simeq\phi-\pi$ near $\phi=\pi$. Thus, the function $T_\nu(z)$ changes rapidly, but smoothly, from being exponentially small in $|\arg\,z|<\pi$ to having the approximate value unity as $\arg\,z$ passes continuously through $\pi$.

The result in (\ref{e35}) and (\ref{e37c}), when the truncation indices $N_k$ and $N_k'$ are chosen according to (\ref{e36}), then constitutes the exponentially improved expansion of $L(\lambda,a,s)$.
For fixed, large $|a|$ in the vicinity of $\arg\,a=\fs\pi$, the dominant contribution to the remainder arises from the term involving $T_{\nu_k'}(iX')$, the other remainder involving $T_{\nu_k}(-iX)$ being smaller. The coefficient of each subdominant exponential $\exp\,(2\pi i(k+\lambda')a)$ then has the leading behaviour from (\ref{e37})  given by
\[e^{-\fr\pi is}\,T_{\nu_k'}(iX')\sim e^{-\fr\pi is}\{\fs+\fs \mbox{erf}\,[c(\theta-\fs\pi)\sqrt{\pi(k+\lambda')|a|}\,]\},\]
where $\theta=\arg\,a$ and, near $\theta=\fs\pi$, the quantity $c(\theta-\fs\pi)\simeq \theta-\fs\pi$. In the vicinity of $\arg\,a=-\fs\pi$, the role of the two remainders is reversed and the coefficient of each subdominant exponential $\exp\,(-2\pi i(k+\lambda)a)$ becomes
\begin{eqnarray}
-e^{-\frac{3}{2}\pi is} T_{\nu_k}(-iX)&=&e^{\fr\pi is} (1-T_{\nu_k}(Xe^{\frac{3}{2}\pi i}))\label{e37b}\\
&\simeq& e^{\fr\pi is}\{\fs-\fs\mbox{erf}\,[c(\theta+\fs\pi)\sqrt{\pi(k+\lambda)|a|}\,]\},\nonumber
\end{eqnarray}
where we have made use of the connection formula for $T_\nu(z)$ given by \cite[Eq.~(6.2.45)]{PK}
\bee\label{e37a}
T_\nu(ze^{-\pi i})=e^{2\pi i\nu} \{T_\nu(ze^{\pi i})-1\}.
\ee
 
The approximate functional form of the Stokes multiplier for $L(\lambda,a,s)$ (excluding the factors $e^{\mp\fr\pi is}/\g(s)$) in the vicinity of $\arg\,a=\pm\fs\pi$ is therefore found to be
\bee\label{e38}
\fs\pm\fs\mbox{erf}\,[(\theta\mp\fs\pi)\sqrt{\pi(k+\xi)|a|}\,],\qquad  (k=0, 1, 2, \ldots ),
\ee
respectively, where $\xi=\lambda'$ near $\arg\,a=\fs\pi$ and $\xi=\lambda$ near $\arg\,a=-\fs\pi$. When $\lambda=1$, $\xi\equiv 0$ and the form (\ref{e38}) then applies to the Hurwitz zeta function $\zeta(s,a)$ with $k=1, 2, \ldots\,$; see \cite[Section 3]{P}.
The form (\ref{e38}) describes the birth of each subdominant exponential in the neighbourhood of the positive and negative imaginary axes on the increasingly sharp scale $(\pi(k+\xi)|a|)^{1/2}$.
It is immediately apparent that the transition across the Stokes lines is associated with {\em unequal scales in the upper and lower half-planes},
except when $\lambda=\fs$ (where the function $L(\fs,a,s)$ reduces to the alternating variant of the Hurwitz zeta function). 
\vspace{0.6cm}

\begin{center}
{\bf 4.\ Numerical results}
\end{center}
\setcounter{section}{4}
\setcounter{equation}{0}
\renewcommand{\theequation}{\arabic{section}.\arabic{equation}}
In order to display numerically the smooth appearance of the $n$th subdominant exponential $e^{2\pi i(n+\lambda')a}$
in the vicinity of $\arg\,a=\fs\pi$ (at fixed $|a|$), it is necessary to `peel off' from $Z(\lambda,a,s)$ the larger subdominant exponentials in the remainder terms and all larger terms of the asymptotic series in (\ref{e35}). 
This has been carried out in the expansion in (\ref{e37c}) by means of the rearrangement in (\ref{e37g})

We define the $n$th Stokes multiplier $S_n(\theta)$ (with $\theta=\arg\,a$) associated with the exponential $e^{2\pi i(n+\lambda')a}$ in the vicinity of $\arg\,a=\fs\pi$ by subtracting from $Z(\lambda,s,a)/(2\pi)^s$ 
the asymptotic series $H_m(a;\lambda,\lambda')$ corresponding to $0\leq m\leq n$ and the
larger subdominant exponentials $0\leq m\leq n-1$ in $R_m'(\lambda',a;N_m')$ and $0\leq m\leq n$ in $R_m(\lambda,a;N_m)$, namely
\[\frac{Z(\lambda,s,a)}{(2\pi)^s}-\frac{i}{2\pi}\sum_{m=0}^n H_m(a;\lambda,\lambda')
=e^{\fr\pi is}\sum_{m=0}^{n-1} \frac{R_m'(\lambda',a;N_m')}{(m+\lambda')^{1-s}}\]
\[-e^{-\fr\pi is}\sum_{m=0}^{n} \frac{R_m(\lambda,a;N_m)}{(m+\lambda)^{1-s}}+\frac{e^{2\pi i(n+\lambda')a-\fr\pi is}}{(n+\lambda')^{1-s}}\,S_n(\theta).
\]
It then follows that near $\arg\,a=\fs\pi$
\[S_n(\theta)=\frac{e^{-2\pi i(n+\lambda')a+\fr\pi is}}{(n+\lambda')^{1-s}}\bl\{\frac{Z(\lambda,s,a)}{(2\pi)^s}-\frac{i}{2\pi}\sum_{m=0}^n H_m(a;\lambda,\lambda')\hspace{4cm}\]
\bee\label{e41}
\hspace{2cm}-e^{\fr\pi is}\sum_{m=0}^{n-1} \frac{R_m'(\lambda',a;N_m')}{(m+\lambda')^{1-s}}
+e^{-\fr\pi is}\sum_{m=0}^{n} \frac{R_m(\lambda,a;N_m)}{(m+\lambda)^{1-s}}\br\}
\ee
for $n=0, 1, 2, \ldots\,$.

Similarly, near $\arg\,a=-\fs\pi$, we write the term corresponding to $m=n$ in the sum involving $R_m(\lambda,a;N_m)$ in (\ref{e37c}) with the aid of (\ref{e37b}). Then the Stokes multiplier near $\arg\,a=-\fs\pi$ is defined by
\[S_n(\theta)=\frac{e^{2\pi i(n+\lambda)a-\fr\pi is}}{(n+\lambda)^{1-s}}\bl\{\frac{Z(\lambda,s,a)}{(2\pi)^s}-\frac{i}{2\pi}\sum_{m=0}^n H_m(a;\lambda,\lambda')\hspace{4cm}\]
\bee\label{e42}
\hspace{2cm}-e^{\fr\pi is}\sum_{m=0}^{n} \frac{R_m'(\lambda',a;N_m')}{(m+\lambda')^{1-s}}
+e^{-\fr\pi is}\sum_{m=0}^{n-1} \frac{R_m(\lambda,a;N_m)}{(m+\lambda)^{1-s}}\br\}
\ee
for $n=0, 1, 2, \ldots\,$.

\begin{table}[t]
\caption{\footnotesize{The real part of the Stokes multiplier $S_0(\theta)$ for $a=5e^{i\theta}$ when $s=4$ and $\lambda=2/3$ compared with the approximate value (\ref{e38}) with $k=0$. The optimal truncation indices are $N_0=17$, $N_0'=7$. }}
\begin{center}
\begin{tabular}{|c|cc||c|cc|}
\hline
&&&&&\\[-0.3cm]
\mcol{1}{|c|}{$\theta/\pi$} & \mcol{1}{c}{$\Re (S_0)$} & \mcol{1}{c||}{Approx $S_0$} & \mcol{1}{c}{$\theta/\pi$} & \mcol{1}{|c}{$\Re (S_0)$} & \mcol{1}{c|}{Approx $S_0$} \\
[.1cm]\hline
&&&&&\\[-0.3cm]
0.30 & 0.02114 & 0.02101 & $-0.30$ & 0.00216 & 0.00202 \\
0.40 & 0.15648 & 0.15466 & $-0.40$ & 0.07660 & 0.07525 \\
0.45 & 0.30653 & 0.30562 & $-0.45$ & 0.23102 & 0.23611 \\
0.48 & 0.41977 & 0.41944 & $-0.48$ & 0.38280 & 0.38685 \\
0.49 & 0.45968 & 0.45951 & $-0.49$ & 0.44063 & 0.44284 \\
0.50 & 0.50000 & 0.50000 & $-0.50$ & 0.50000 & 0.50000 \\
0.51 & 0.54032 & 0.54049 & $-0.51$ & 0.55937 & 0.55716 \\
0.52 & 0.58023 & 0.58056 & $-0.52$ & 0.61720 & 0.61315 \\
0.55 & 0.69347 & 0.69438 & $-0.55$ & 0.76898 & 0.76389 \\
0.60 & 0.84352 & 0.84534 & $-0.60$ & 0.92340 & 0.92475 \\
0.70 & 0.97886 & 0.97899 & $-0.70$ & 0.99784 & 0.99798 \\
[.2cm]\hline
\end{tabular}
\end{center}
\end{table}
\begin{table}[t]
\caption{\footnotesize{The real part of the Stokes multiplier $S_1(\theta)$ for $a=5e^{i\theta}$ when $s=4$ and $\lambda=2/3$ compared with the approximate value (\ref{e38}) with $k=1$. The optimal truncation indices are $N_0=17$, $N_0'=7$, $N_1=49$, $N_1'=38$. }}
\begin{center}
\begin{tabular}{|c|cc||c|cc|}
\hline
&&&&&\\[-0.3cm]
\mcol{1}{|c|}{$\theta/\pi$} & \mcol{1}{c}{$\Re (S_1)$} & \mcol{1}{c||}{Approx $S_1$} & \mcol{1}{c}{$\theta/\pi$} & \mcol{1}{|c}{$\Re (S_1)$} & \mcol{1}{c|}{Approx $S_1$} \\
[.1cm]\hline
&&&&&\\[-0.3cm]
0.35 & 0.00114 & 0.00114 & $-0.35$ & 0.00021 & 0.00032 \\
0.40 & 0.02157 & 0.02101 & $-0.40$ & 0.01128 & 0.01151 \\
0.45 & 0.15510 & 0.15466 & $-0.45$ & 0.12807 & 0.12785 \\
0.48 & 0.34208 & 0.34213 & $-0.48$ & 0.32480 & 0.32468 \\
0.49 & 0.41939 & 0.41944 & $-0.49$ & 0.41014 & 0.41009 \\
0.50 & 0.50000 & 0.50000 & $-0.50$ & 0.50000 & 0.50000 \\
0.51 & 0.58061 & 0.58056 & $-0.51$ & 0.58986 & 0.58991 \\
0.52 & 0.65792 & 0.65787 & $-0.52$ & 0.67520 & 0.67532 \\
0.55 & 0.84490 & 0.84534 & $-0.55$ & 0.87193 & 0.87215 \\
0.60 & 0.97843 & 0.97899 & $-0.60$ & 0.98872 & 0.98849 \\
0.65 & 0.99889 & 0.99886 & $-0.65$ & 0.99979 & 0.99968 \\
[.2cm]\hline
\end{tabular}
\end{center}
\end{table}

In Tables~1 and 2 we show the real part\footnote{There is also a small imaginary part to $S_n(\theta)$ that is not presented.} of $S_n(\theta)$ for $n=0$ and $1$ computed from (\ref{e41}) and (\ref{e42}) compared with the  approximate value in (\ref{e38}) when $a=5e^{i\theta}$ and $s=4$, $\lambda=\f{2}{3}$ as a function of $\theta$ in the vicinity of the positive and negative imaginary $a$-axes. 
In the computation of $S_n(\theta)$ it is necessary to compute the terms $R_m(\lambda,a;N_m)$ and $R_m'(\lambda',a;N_m')$ by means of the incomplete gamma function representation in (\ref{e32a}) and the sum $H_m(a;\lambda,\lambda')$ to the required exponential accuracy. The optimal truncation indices
$N_k$ and $N'_k$ were obtained by inspection of the terms in the algebraic expansions.
In addition, when computing the terminant functions appearing in $R_m(\lambda,a;N_m)$ and $R_m'(\lambda',a;N_m')$ one must use the connection formula (\ref{e37b}) once the argument of $z$ in $T_\nu(z)$ has exceeded $\pi$, since {\it Mathematica} only computes the value of the incomplete gamma function in the principal sector $-\pi<\arg\,z\leq\pi$.

It is seen that there is good agreement between the real part of the computed values of the Stokes multiplier and the predicted approximate values in (\ref{e38}). Moreover, the tables confirm that the transition scales across the positive and negative imaginary $a$-axes depend on $\lambda'$ and $\lambda$, respectively, and are indeed unequal (except when $\lambda=\fs$).

\vspace{0.6cm}

\begin{center}
{\bf Appendix A: \ Estimation of the remainder ${\cal R}_K(\lambda,a,s)$ in (\ref{e24})}
\end{center}
\setcounter{section}{1}
\setcounter{equation}{0}
\renewcommand{\theequation}{\Alph{section}.\arabic{equation}}
The remainder term ${\cal R}_K(\lambda,a,s)$ in (\ref{e24}) resulting from displacement of the integration path in (\ref{e23}) is given by
\[{\cal R}_K(\lambda,a,s)=\frac{a^{-s}}{2\pi i}\int_{K-c-\infty i}^{K-c+\infty i} \g(-u) \g(s+u) F(\lambda,-u)\,a^{-u} du\qquad (0<c<1).\]
Following the procedure described in Section 3 combined with use of the functional equation for $F(\lambda,-u)$ in (\ref{e30}), we easily obtain
\bee\label{a1}
\frac{{\cal R}_K(\lambda,a,s)}{(2\pi)^s}=e^{-\fr\pi is}\sum_{k=0}^\infty{}{\!'}(k+\lambda')^{s-1} e^{iX'}\,T_\nu(iX')-e^{-\frac{3}{2}\pi is}
\sum_{k=0}^\infty(k+\lambda)^{s-1} e^{-iX}\,T_\nu(-iX),
\ee
where $T_\nu(z)$ is the terminant function in (\ref{e32a}), $\nu=K+s$ and $X$, $X'$ are defined in (\ref{e32d}).

Since \cite[p.~260]{PK}
\[e^zT_\nu(z)=-ie^{\pi i\nu}\,\frac{\g(\nu)}{2\pi}\,U(\nu,\nu,z),\]
where $U(a,b,z)$ denotes the second confluent hypergeometric function \cite[p.~322]{DLMF}, it is seen that the expansion (\ref{e24}) with the remainder given in (\ref{a1}) becomes
\[L(\lambda,a,s)=a^{-s}+\sum_{k=0}^{K-1}\frac{(-)^k}{k!} (s)_k F(\lambda,-k)\,a^{-s-k}+(-)^K(2\pi)^s
\,\frac{(s)_K}{2\pi i}\hspace{2cm}\]
\bee\label{a2}
\hspace{1cm}\times\bl\{e^{\fr\pi is}\sum_{k=0}^\infty{}(k+\lambda')^{s-1} U(\nu,\nu,iX')-e^{-\fr\pi is}\sum_{k=0}^\infty (k+\lambda)^{s-1}U(\nu,\nu,-iX)\br\}
\ee
when $0<\lambda<1$.
This is similar to, but not identical with, the expansion obtained in \cite{K}.

For fixed $\nu$, we have the asymptotic behaviour \cite[p.~328]{DLMF}
\[U(\nu,\nu,z)\sim z^{-\nu} \qquad (|z|\ra\infty,\ |\arg\,z|<\f{3}{2}\pi).\]
It then follows that for {\it fixed\/} integer $K$
\begin{eqnarray}
\frac{{\cal R}_K(\lambda,a,s)}{(2\pi)^s}&\sim&(2\pi a)^{-\nu}\,\frac{\g(\nu)}{2\pi i}\bl\{e^{\fr\pi iK} \sum_{k=0}^\infty{}\!' (k+\lambda')^{-K-1}-e^{-\fr\pi iK}\sum_{k=0}^\infty (k+\lambda)^{-K-1}\br\}\nonumber\\
&=&O(a^{-K-s}) \label{a3}
\end{eqnarray}
 as $|a|\ra\infty$ in $|\arg\,a|<\pi$, since the sums are finite and independent of $a$.
\vspace{0.6cm}

\begin{center}
{\bf Appendix B: \ The double series rearrangement in (\ref{e37g})}
\end{center}
\setcounter{section}{2}
\setcounter{equation}{0}
\renewcommand{\theequation}{\Alph{section}.\arabic{equation}}
The rearrangement of the double series in (\ref{e37g}) follows the procedure described in \cite{B2}; see also \cite[\S 6.4.3]{PK} for an account of this process. 
If we set 
$A_r:=(-i)^r \g(r+s)/(2\pi a)^{r+s}$,
the double sum in (\ref{e35}) involving $\lambda$ and the truncation indices $\{N_k\}$ ($k\geq 0$) can be rearranged\footnote{The modified double series involves the function $\zeta(r+1,m+\lambda)$; however, its evaluation is straightforward for $m=0, 1, 2, \ldots\,$.} as
\[\sum_{k=0}^\infty \sum_{r=1}^{N_k-1} \frac{A_r}{(k+\lambda)^{r+1}}\hspace{11cm}\]
\begin{eqnarray*}
&=&\sum_{r=1}^{N_0-1}\frac{A_r}{\lambda^{r+1}}+\bl(\sum_{r=1}^{N_0-1}+\sum_{r=N_0}^{N_1-1}\br)\frac{A_r}{(1+\lambda)^{r+1}}+
\bl(\sum_{r=1}^{N_0-1}+\sum_{r=N_0}^{N_1-1}+\sum_{r=N_1}^{N_2-1}\br)\frac{A_r}{(2+\lambda)^{r+1}}+ \cdots\\
&=&\sum_{r=1}^{N_0-1}A_r\sum_{k=0}^\infty \frac{1}{(k+\lambda)^{r+1}}+\sum_{r=N_0}^{N_1-1}A_r\sum_{k=1}^\infty \frac{1}{(k+\lambda)^{r+1}}
+\sum_{r=N_1}^{N_2-1}A_r\sum_{k=2}^\infty \frac{1}{(k+\lambda)^{r+1}}+\cdots\\
&=&\sum_{m=0}^\infty \sum_{r=N_{m-1}}^{N_m-1} A_r\,\zeta(r+1,m+\lambda),
\end{eqnarray*}
where the sums over $k$ have been evaluated in terms of the Hurwitz zeta function and $N_{-1}=1$. A similar result applies to the other double sum involving $\lambda'$ in (\ref{e35}) with $N_{-1}'=1$.

\end{document}